\documentclass[12pt]{article}

\usepackage{amsmath, amssymb, amsfonts, mathrsfs, enumitem, tikz-cd, mathtools}
\usepackage{comment}
\usepackage{graphicx}
\usepackage{xcolor}
\usepackage{tikz}
\usepackage[framemethod=TikZ]{mdframed}
\usetikzlibrary{backgrounds}
\usepackage{longtable}
\usepackage[most]{tcolorbox}
\usepackage[colorlinks=true,linkcolor=black,anchorcolor=black,citecolor=black,filecolor=black,menucolor=black,runcolor=black,urlcolor=blue]{hyperref}
\usetikzlibrary{external, fit, shapes, decorations.pathreplacing, calligraphy, arrows.meta, positioning}
\usetikzlibrary{decorations.markings,arrows}

\usepackage[margin=1.13 in, top=1.1 in, bottom= 1.2 in]{geometry}
\newcommand{\N}{\mathbb{N}}

\newcommand{\C}{\mathbb{C}}
\newcommand{\oo}{\infty}
\newcommand{\sub}{\subseteq}

\usepackage{amsthm}
\newtheoremstyle{plain}{3mm}{3mm}{\slshape}{}{\bfseries}{.}{.5em}{}
\newtheoremstyle{definition}{2mm}{2mm}{}{}{\bfseries}{.}{.5em}{}
\theoremstyle{plain} 
	
\newtheorem{theorem}{Theorem}[section]

\newtheorem{lemma}[theorem]{Lemma}
\newtheorem{corollary}[theorem]{Corollary}
\theoremstyle{definition} 
\newtheorem{definition}[theorem]{Definition}

\theoremstyle{plain} 
\newcounter{MainTheoremCounter}

\theoremstyle{plain}
\newtheorem*{namedthm}{\namedthmname}
\newcounter{namedthm}
\makeatletter
	\newenvironment{named}[2]
	{\def\namedthmname{#1}
	\refstepcounter{namedthm}
	\namedthm[#2]\def\@currentlabel{#1}}
	{\endnamedthm}
\makeatother

\title{Recurrence Relations for Cosets in Free Groups}

\author{Michael Reilly and Cory Shields}

\begin{document}

\maketitle
\begin{abstract}
    Let $F_2$ be the free group on two generators and let $H$ be a subgroup of $F_2$. We investigate a method for calculating the number of elements in a coset of $H$ that have a given length when written in reduced form. More specifically, taking $S_n\subseteq F_2$ to be the set of elements of length $n$, we show that for any coset $yH$ there always exists a recurrence relation of the form
    \begin{equation*}
          |yH\cap S_n| = \sum_{i=1}^{n-1}\sum_{xH\in F_2/H}a_{i,xH}\cdot |xH\cap S_{n-i}|
    \end{equation*}
    for some constants $(a_{i,xH})_{i\in \mathbb{N}, xH\in F_2/H}$, and we give an algorithm that calculates these constants. Further, we show that when $H$ has finite index and contains an element of odd length, only finitely many of the constants $a_{i,xH}$ are nonzero.
\end{abstract}

\section{Introduction}\label{section:introduction}

Let $F_2$ denote the free group on two generators, $a$ and $b$. Let $\Xi = \{a,b,a^{-1},b^{-1}\}$. We say that $g\in F_2$ has \emph{length} $n$ if $g = x_1x_2\dots x_n$ where $x_i\in \Xi$ and $x_{i+1}\neq x_i^{-1}$ for each $i$. For each $n\geq 1$, let $S_n\subseteq F_2$ be the set of all elements of length $n$. Our goal is to be able to calculate $|yH\cap S_n|$, where $H$ is a subgroup of $F_2$ and $y\in F_2$. We do this by finding a recurrence relation for $|yH\cap S_n|$ in terms of $|xH\cap S_k|$ for $xH\in F_2/H$ and $k\leq n-1$.

This leads us to our main theorems, whose proofs are given in Sections \ref{section:proof_of_main_thm} and \ref{section:rec_relation_is_finite}, respectively.
\begin{theorem}\label{thm:main}
Let $H\leq  F_2$ and let $y\in F_2$. There exist constants $a_{i,xH}$ for $i\geq 1$ and $xH\in F_2/H$ such that for each $n\geq 2$,
\begin{align*}
    |yH\cap S_n| = \sum_{i=1}^{n-1}\sum_{xH\in F_2/H}a_{i,xH}\cdot |xH\cap S_{n-i}|.
\end{align*}
\end{theorem}

\begin{theorem}\label{thm:main_2}

Let $H\leq  F_2$ and let $y\in F_2$. Further, suppose $H$ has finite index and $H$ contains an element of odd length. Then there is an $N\in \N$ such that 
\begin{align*}
    |yH\cap S_n| = \sum_{i=1}^{N}\sum_{xH\in F_2/H}a_{i,xH}\cdot |xH\cap S_{n-i}|
    \end{align*}
    for all $n\geq N+1$.
\end{theorem}

We also give an algorithm for computing the coefficients $(a_{i,xH})$ in our main theorems. Given a subgroup $H\leq F_2$, the \emph{Directed Coset Graph} of $F_2/H$ is a directed graph with a vertex set equal to $F_2/H$ and with a directed edge from $rH$ to $sH$ for each $h\in \Xi$ such that $hrH = sH$. A directed coset graph could have edges from a vertex to itself or two or more edges between two vertices.

To compute the coefficients $(a_{i,xH})$, perform the following algorithm:

Fix a subgroup $H\leq F_2$ and an element $y\in F_2$. We will say that a directed edge is \emph{highlighted} in a given step to mean that it belongs to the set of edges that we will use to calculate the coefficients in that step. We will say that a vertex $xH$ is \emph{charged on Step $i$} to mean that in Step $i$ there is a highlighted directed edge that points into $xH$. By convention, we will say that $yH$ is charged on Step $0$.

\begin{enumerate}
    \item In Step 1 highlight all directed edges that point out of the vertex $yH$.
    \item At the end of Step $i\geq 1$, for each $xh\in F_2/H$ count how many highlighted directed edges point into $xH$ and call this number $b_{i,xH}$. Put $a_{i,xH} = b_{i,xH}-1$ if $b_{i,xH}>0$ and $a_{i,xH} = 0$ otherwise. Then unhighlight all highlighted directed edges.
    \item At the start of Step $(i+1)$ with $i\geq 1$, for each vertex $xH\in F_2/H$ which was charged on Step $i$, highlight each directed edge pointing out of $xH$ except for edges that point into a vertex that was charged on Step $i-1$.
\end{enumerate}

Note that an edge $rH\to sH$ is highlighted on Step $i+1$ iff $rH$ is charged on Step $i$ and $sH\rightarrow rH$ is not highlighted on Step $i$. As one might expect, the presence of parallel edges, which are distinct edges that both begin and end at the same pair of vertices, are counted in the coefficients $a_{i,xH}$ but they do not affect which cosets are charged. Presently, it is unclear whether this algorithm will ever terminate. In Section \ref{section:rec_relation_is_finite}, we will see that when $H$ is a finite index subgroup that contains an element of odd length, this algorithm always ends after finitely many steps.

Examples of this algorithm are given in Section \ref{section:examples}.

\subsection*{Acknowledgments}
This work was done as a part of the Cycle undergraduate research program at The Ohio State University. The authors would like to thank this program for its generous nonfinancial support, without which this work would not have been possible.

\section{Examples}\label{section:examples}

Let us see some examples of using the algorithm described above to calculate recurrence relations. In the examples that follow, we will make highlighted edges green, and for the purpose of having less visual clutter, we will neglect to draw the arrowhead on edges that are not highlighted. When we introduce the coset graph for each subgroup, we will label each arrowhead on a directed edge with the element of $\Xi$ that it represents. For the purpose of avoiding visual clutter once more on the coset graph, we choose to place two arrowheads labeled $g$ and $g^{-1}$ respectively for some $g\in\Xi$  on the same edge instead of two separate edges. Additionally, we will keep track of the coefficients $(a_{i,xH})$ by writing $a_{i,xH}$ as a superscript of $xH$ in Step $i$ whenever $a_{i,xH}\neq 0$. 

\textbf{Example 1}:

  Let $\pi:F_2\rightarrow S_3$ be defined by $\pi(a) = (12)$ and $\pi(b) = (23)$, and let $H=\ker \pi$. $H$ can also be written as $H = \langle a^2, b^2, a^{-1}ba, b^{-1}ab\rangle$.
Then $H$ has the following coset graph.

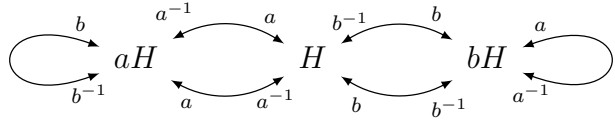
\begin{center}
\begin{tikzcd}[row sep=1cm, column sep=1.5cm]
   \arrow[loop left, distance = 4em, no head, in = 160, out = 200, "b^{-1}"{inner sep=1pt, pos=0.1}, "b" pos=0.9,latex-latex] aH & \arrow[l,no head, bend left, "a^{-1}"{inner sep=-.1pt, near start}, "a" near end,latex-latex] \arrow[l,no head, bend right, "a"{ near start}, "a^{-1}" near end,swap, latex-latex] H \arrow[r,no head, bend left, "b^{-1}"{inner sep=-.1pt, near start}, "b" near end,latex-latex] \arrow[r,no head, bend right, "b"{ near start}, "b^{-1}" near end,swap, latex-latex] & bH   \arrow[loop right, distance = 4em, no head, in = 20, out = -20, "a^{-1}"{inner sep=1pt, pos=0.1}, "a" pos = 0.9,latex-latex, swap] 
\end{tikzcd}
\end{center}
We can compute the recurrence relation for $|aH\cap S_n|$ using these steps below.

\begin{center}

\begin{tikzcd}
    \draw node [draw, ultra thick] {\textbf{Step 1}};  & & 
    \\
   \arrow[loop left, distance  = 3em, no head, in = 160, out = 200, thick, green, latex-latex] aH^{(1)} & \arrow[l,no head, bend left,  thick, green, latex-] \arrow[l,no head, bend right,  thick, green, latex-] H^{(1)} \arrow[r,no head, bend left] \arrow[r,no head, bend right] & bH   \arrow[loop right, distance  = 3em, no head] 
\end{tikzcd}

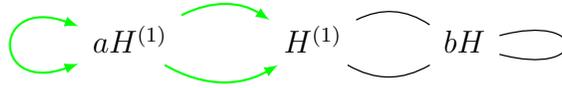
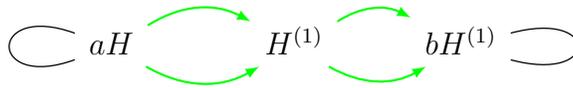
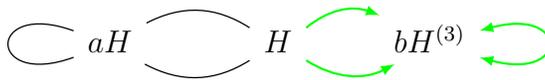
\begin{tikzcd}
    \draw node [draw, ultra thick] {\textbf{Step 2}};  & & 
    \\
   \arrow[loop left, distance  = 3em, no head, in = 160, out = 200] aH & \arrow[l,no head, bend left, thick, green, latex-] \arrow[l,no head, bend right, thick, green, latex-] H^{(1)} \arrow[r,no head, bend left, thick, green, -latex] \arrow[r,no head, bend right, thick, green, -latex] & bH^{(1)}   \arrow[loop right, distance  = 3em, no head] 
\end{tikzcd}

\begin{tikzcd}
    \draw node [draw, ultra thick] {\textbf{Step 3}};  & & 
    \\
   \arrow[loop left, distance  = 3em, no head, in = 160, out = 200,] aH & \arrow[l,no head, bend left] \arrow[l,no head, bend right] H \arrow[r,no head, bend left, thick, green, -latex] \arrow[r,no head, bend right, thick, green, -latex] & bH^{(3)}  \arrow[loop right, distance  = 3em, no head, thick, green, latex-latex] 
\end{tikzcd}

\end{center}

From this we calculate that, 
\begin{align*}
    |aH\cap S_n| = |aH\cap S_{n-1}|+|H\cap S_{n-1}|+|H\cap S_{n-2}|+|bH\cap S_{n-2}|+3|bH\cap S_{n-3}|
\end{align*}
for each $n\geq 3$.

\textbf{Example 2}:

Let $\pi: F_2\to S_{5}$ be defined by $\pi(a) = (12)(45)$ and $\pi(b) = (14)(235)$, and let $H=\ker \pi$. $H$ can also be written as $H=\langle a^2, b^2, (ba)^{-1}aba, (ab)^{-1}b^2a,b^{-1}a^2b,a^{-1}bab \rangle$. Then $H$ has the following coset graph.

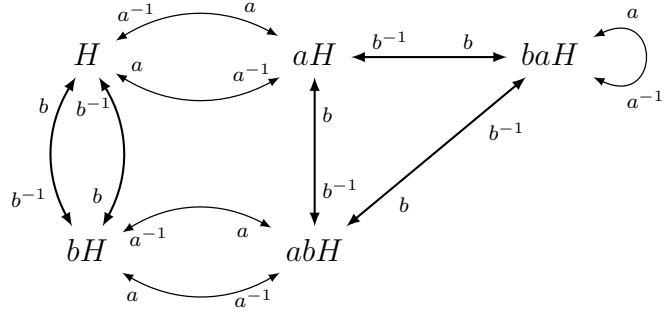
\begin{center}

\begin{tikzcd}[row sep=2cm, column sep=2cm]
    H\arrow[r, no head, bend left, "a^{-1}"{inner sep=-.1pt, near start}, "a" near end,latex-latex]\arrow[r, no head, bend right, "a^{-1}"{pos=.95, inner sep=-.1pt}, "a"{pos=.05}, latex-latex]\arrow[d, no head, bend left, thick, "b^{-1}"{pos=.15,  inner sep=-.9pt, swap}, "b" {pos=.9,  swap}, latex-latex]\arrow[d, no head, bend right, thick, "b^{-1}"{pos=.7, swap}, "b" {pos=.295,  swap}, latex-latex]
    & aH\arrow[r, no head, thick, "b^{-1}" near start, "b" near end, latex-latex]\arrow[d, no head, thick, "b" near start, "b^{-1}" near end, latex-latex]
    & baH\arrow[loop above, out=25, in=-25, distance=1cm,
               no head,
               "a_{\phantom{-1}}" {pos = 0.15}, "a^{-1}" {pos = 0.85}, latex-latex]\arrow[dl, no head, thick, "b^{-1}" near start, "b" near end, latex-latex]\\
    bH\arrow[r, no head, bend left, "a^{-1}"{pos=.05, inner sep=-.3pt, swap}, "a"{pos=.9, swap},latex-latex]\arrow[r, no head, bend right, "a^{-1}"{pos=.725, swap, inner sep=-.9pt}, "a"{pos=.15, swap}, latex-latex]
    & abH &
\end{tikzcd}
\end{center}

We can compute the recurrence relation for $|H\cap S_n|$ using the following steps.

\begin{tikzcd}
    \draw node [draw, ultra thick] {\textbf{Step 1}};  & & 
    \\
    H\arrow[r, no head, bend right=20, thick, green, -latex]\arrow[r, no head, bend left=27,thick, green, -latex]\arrow[d, no head, bend right=25, thick, green, -latex]\arrow[d, no head, bend left=25,thick, green, -latex] & aH^{(1)}\arrow[r, no head]\arrow[d, no head] & 
    baH\arrow[loop right, distance = 2.5em, no head, in = 30, out = -30]\arrow[dl, no head]\\
    bH^{(1)}\arrow[r, no head, bend right]\arrow[r, no head, bend left] & abH &
\end{tikzcd}
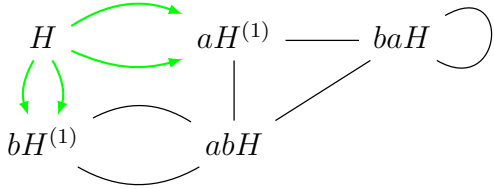
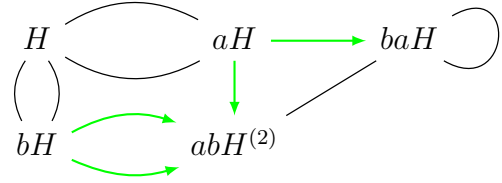
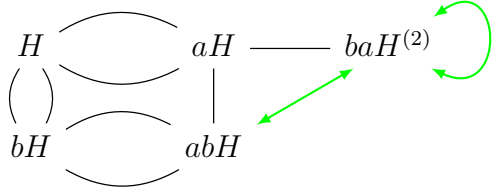
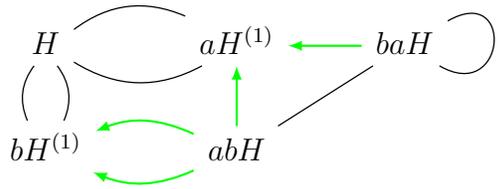
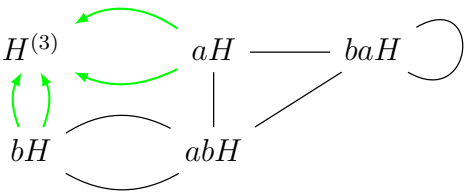
\begin{tikzcd}
\draw node [draw, ultra thick] {\textbf{Step 2}};  & & 
    \\
    H\arrow[r, no head, bend right]\arrow[r, no head, bend left]\arrow[d, no head, bend right]\arrow[d, no head, bend left] & aH\arrow[r, no head, thick, green, -latex]\arrow[d, no head, thick, green, -latex] & 
    baH\arrow[loop right, distance = 2.5em, no head, in = 30, out = -30]\arrow[dl, no head]\\
    bH\arrow[r, no head, bend right=20, thick, green, -latex]\arrow[r, no head, bend left=23,thick, green, -latex] & abH^{(2)} &
\end{tikzcd}

\hfill \phantom{.}

\begin{tikzcd}
\draw node [draw, ultra thick] {\textbf{Step 3}};  & & 
    \\
    H\arrow[r, no head, bend right]\arrow[r, no head, bend left]\arrow[d, no head, bend right]\arrow[d, no head, bend left] & aH\arrow[r, no head]\arrow[d, no head] & baH^{(2)}\arrow[loop right, distance = 2.75em, no head, in = 35, out = -25, thick, green, latex-latex]\arrow[dl, no head, thick, green, latex-latex]\\
    bH\arrow[r, no head, bend right]\arrow[r, no head, bend left] & abH &
\end{tikzcd}
\begin{tikzcd}
\draw node [draw, ultra thick] {\textbf{Step 4}};  & & 
    \\
    H\arrow[r, no head, bend right]\arrow[r, no head, bend left]\arrow[d, no head, bend right]\arrow[d, no head, bend left] & aH^{(1)}\arrow[r, no head, thick, green, latex-]\arrow[d, no head, thick, green, latex-] & baH\arrow[loop right, distance = 2.5em, no head, in = 30, out = -30]\arrow[dl, no head]\\
    bH^{(1)}\arrow[r, no head, bend right=23, thick, green, latex-]\arrow[r, no head, bend left=23, thick, green, latex-] & abH &
\end{tikzcd}
\hfill \phantom{.}

\begin{tikzcd}
\draw node [draw, ultra thick] {\textbf{Step 5}};  & & 
    \\
    H^{(3)}\arrow[r, no head, bend right=25, thick, green, latex-]\arrow[r, no head, bend left=33,thick, green, latex-]\arrow[d, no head, bend right=25,thick, green, latex-]\arrow[d, no head, bend left=25,thick, green, latex-] & aH\arrow[r, no head]\arrow[d, no head] & baH\arrow[loop right, distance = 2.5em, no head, in = 30, out = -30]\arrow[dl, no head]\\
    bH\arrow[r, no head, bend right]\arrow[r, no head, bend left] & abH &
\end{tikzcd}

As a result,
\begin{align*}
    |H\cap {S}_{n}|=~&|aH\cap S_{n-1}|+|bH\cap S_{n-1}|+2|abH\cap S_{n-2}|+2|baH\cap S_{n-3}|\\
    +&|aH\cap S_{n-4}|+|bH\cap S_{n-4}|+3|H\cap S_{n-5}|\quad \text{ for each } n\geq 5.
    \end{align*}

\textbf{Example 3:} 
Let $\pi:F_2\rightarrow S_7$ be defined by $\pi(a) = (215)(47)$ and $\pi(b) = (15)(23476)$, and let $H=\ker \pi$. $H$ can also be written as $H=\langle a^3, b^5, aba, bab^{-1},b^{-1}ab,a^{-1}ba^{-1},b^2ab^2,b^{-2}ab^{-2}\rangle$.
Then $H$ has the following coset graph.

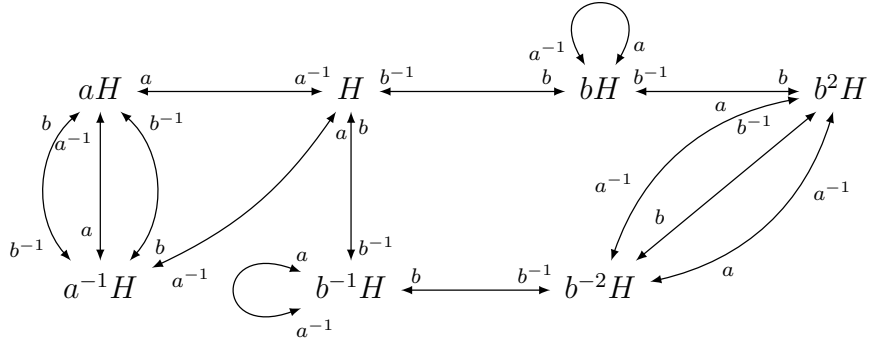
\begin{center}
\begin{tikzcd}[row sep=2cm, column sep=2cm]
aH \arrow[d, latex-latex, "a^{-1}"{pos = 0.2}, "a"pos = 0.8,swap]\arrow[d, bend left = 45, latex-latex,"b^{-1}"{pos = 0.2}, "b"{pos = 0.8}]\arrow[d, bend right = 45, latex-latex,"b"{pos = 0.2}, "b^{-1}"{pos = 0.8},swap]& H\arrow[l, latex-latex , "a^{-1}"{pos = 0.05}, "a"pos = 0.95,swap]\arrow[dl, latex-latex,"a"{pos = 0.05}, "a^{-1}" pos = 0.95, bend left = 15]\arrow[d, latex-latex,"b"{pos = 0.1}, "b^{-1}"pos = 0.9]\arrow[r, latex-latex,"b^{-1}"{pos = 0.1}, "b"pos = 0.9] & bH\arrow[loop, in = 60, out = 120, distance = 3em, latex-latex, "a^{-1}"{pos = 0.15}, "a"pos = 0.85] \arrow[r, latex-latex,"b^{-1}"{pos = 0.1}, "b"pos = 0.9]& b^2H\arrow[dl, latex-latex, "b^{-1}"{pos = 0.2}, "b"pos = 0.8,swap]\arrow[dl, latex-latex, bend left=30,"a^{-1}"{pos = 0.25}, "a"pos = 0.75]\arrow[dl, latex-latex, bend right= 30,"a"{pos = 0.25}, "a^{-1}"pos = 0.75,swap ] \\
    a^{-1}H & \arrow[loop, in = 160, out = 200, distance = 3em, latex-latex,"a^{-1}"{pos = 0.05}, "a"pos = 0.95]b^{-1}H \arrow[r, latex-latex,"b"{pos = 0.1}, "b^{-1}"pos = 0.9]& b^{-2}H
\end{tikzcd}
\end{center}
We can compute the recurrence relation for $|H\cap S_n|$ using the following steps.

\begin{center}
\begin{tikzcd}
\draw node [draw, ultra thick] {\textbf{Step 1}};  & & & \\
    aH \arrow[d, no head]\arrow[d, bend left, no head]\arrow[d, bend right, no head]& H\arrow[l, no head, thick, green, -latex]\arrow[dl, no head,  thick, green, -latex]\arrow[d, no head,  thick, green, -latex]\arrow[r, no head,  thick, green, -latex] & bH\arrow[loop, in = 60, out = 120, distance = 3em, no head] \arrow[r, no head]& b^2H\arrow[dl, no head]\arrow[dl, no head, bend left=15]\arrow[dl, no head, bend right= 15] \\
    a^{-1}H & \arrow[loop, in = 160, out = 200, distance = 3em, no head]b^{-1}H \arrow[r, no head]& b^{-2}H
\end{tikzcd}

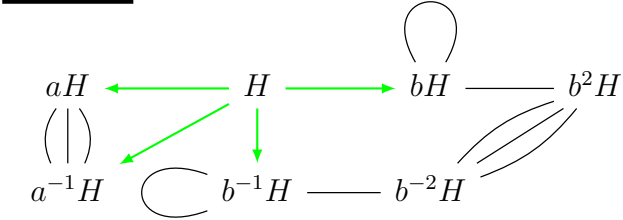
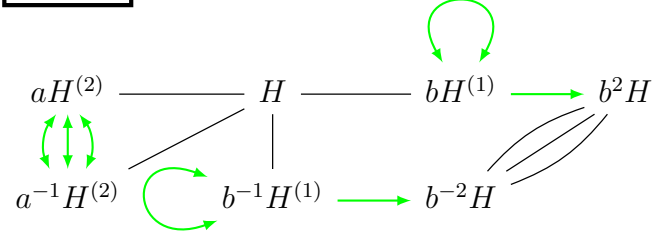
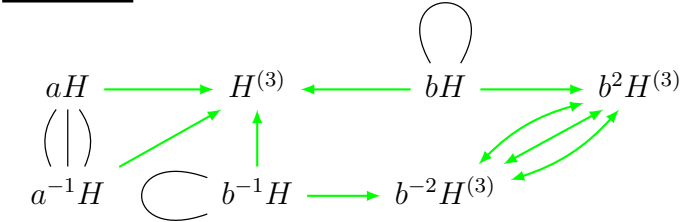
\begin{tikzcd}
\draw node [draw, ultra thick] {\textbf{Step 2}};  & & & \\
    aH^{(2)} \arrow[d, no head,thick, green, latex-latex]\arrow[d, bend left, no head, thick, green, latex-latex]\arrow[d, bend right, no head, thick, green, latex-latex]& H\arrow[l, no head]\arrow[dl, no head]\arrow[d, no head]\arrow[r, no head] & bH^{(1)}\arrow[loop, in = 60, out = 120, distance = 3em, no head,thick, green, latex-latex] \arrow[r, no head,thick, green, -latex]& b^2H\arrow[dl, no head]\arrow[dl, no head, bend left=15]\arrow[dl, no head, bend right= 15]  \\
    a^{-1}H^{(2)} &\arrow[loop left, no head, in = 160, out = 200, distance = 3em,thick, green, latex-latex]b^{-1}H^{(1)} \arrow[r, no head,thick, green, -latex]& b^{-2}H
\end{tikzcd}

\begin{tikzcd}
\draw node [draw, ultra thick] {\textbf{Step 3}};  & & & \\
    aH \arrow[d, no head]\arrow[d, bend left, no head]\arrow[d, bend right, no head]& H^{(3)}\arrow[l, no head, thick, green, latex-]\arrow[dl, no head, thick, green, latex-]\arrow[d, no head, thick, green, latex-]\arrow[r, no head,thick, green, latex-] & bH\arrow[loop, in = 60, out = 120, distance = 3em, no head] \arrow[r, no head,  thick, green, -latex]& b^2H^{(3)}\arrow[dl, no head,  thick, green, latex-latex]\arrow[dl, no head, bend left=15,  thick, green, latex-latex]\arrow[dl, no head, bend right =15,  thick, green, latex-latex] \\
    a^{-1}H & \arrow[loop left, no head, in = 160, out = 200, distance = 3em]b^{-1}H \arrow[r, no head, thick, green, -latex]& b^{-2}H^{(3)}
\end{tikzcd}
\end{center}

From this we calculate that,
\begin{align*}
    |H\cap S_{n}|=~&2|aH\cap S_{n-2}|+2|a^{-1}H\cap S_{n-2}|+|bH\cap S_{n-2}|+|b^{-1}H\cap S_{n-2}|\\+&3|H\cap S_{n-3}|+3|b^2H\cap S_{n-3}|+3|b^{-2}H\cap S_{n-3}|\quad\text{ for each } n\geq 3.
    \end{align*}

\section{Proof of Theorem \ref{thm:main}} 
\label{section:proof_of_main_thm}

Now we shall set out to prove the main theorem.

\begin{named}{Theorem \ref{thm:main}}{}
Let $H\leq  F_2$ and let $y\in F_2$. There exist constants $a_{i,xH}$ for $i\geq 1$ and $xH\in F_2/H$ such that for each $n\geq 2$, we have
\begin{align*}
    |yH\cap S_n| = \sum_{i=1}^{n-1}\sum_{xH\in F_2/H}a_{i,xH}\cdot |xH\cap S_{n-i}|.
    \end{align*}
\end{named}

Before proceeding, we introduce some helpful facts. Recall that $\Xi = \{a,b,a^{-1},b^{-1}\}$. For each $g\in \Xi$, let $A_{g}\subseteq F_2\setminus \{e\}$ be the set of elements that do not start with $g$ when written in reduced form from left to right. For example, $a^{-3},bab,b^{-17}a\in A_{a}$, but $a,ab^{-10},aba^{-1}\not\in A_a$. Each nonidentity element of $F_2$ is contained in exactly $3$ of $A_a,A_b,A_{a^{-1}},A_{b^{-1}}$. Equivalently, we could consider the complements of these sets, taken with respect to $F_2\setminus \{e\}$, so that $A_a^c,A_b^c,A_{a^{-1}}^c,A_{b^{-1}}^c$ forms a partition of $F_2\setminus \{e\}$.

\begin{lemma}
    Let $H\leq F_2$. Then for any $xH\in F_2/H$, $n\geq 2$ and $h\in \Xi$,
    \begin{equation}\label{eq:also_inductive_step}
    |xH\cap S_n\cap A_{h}^c| = |h^{-1}xH\cap S_{n-1}\cap A_{h^{-1}}|.
    \end{equation}
\end{lemma}

\begin{proof}

    Consider the map $g\mapsto h^{-1}g$, which defines an automorphism of $F_2$. Each element of $xH$ is mapped to an element of $h^{-1}xH$. An element of $ S_n\cap A_{h}^c$ can be written in reduced form as $x_1x_2\dots x_n$, where each $x_i\in \Xi$, $x_1 = h$ and $x_ix_{i+1}^{-1}\neq e$ for all $i$. Hence, $h^{-1}x_1x_2\dots x_n = x_2\dots x_n$, which is an element of length $n-1$ that starts with $x_2\neq h^{-1}$. Therefore, each element of $S_n\cap A_{h}^c$ is mapped to an element of $S_{n-1}\cap A_{h^{-1}}$. Conversely, the same argument shows that each element of $S_{n-1}\cap A_{h^{-1}}$ is mapped to $S_n\cap A_{h}^c$ under the inverse map $g\mapsto hg$. Thus, we can conclude that $|xH\cap S_n\cap A_{h}^c| = |h^{-1}xH\cap S_{n-1}\cap A_{h^{-1}}|$.
    
\end{proof}

\begin{lemma}
Let $H\leq F_2$. Then for any $xH\in F_2/H$ and $n\geq 2$,
    \begin{equation}\label{rec_relation}
        |xH\cap S_n| = \sum_{h\in \Xi}|{hx}H\cap S_{n-1}\cap A_{h}|.
    \end{equation}
\end{lemma}

\begin{proof}

Since $\bigcup_{h\in \Xi}A_h^c$ forms a partition of $F_2\setminus \{e\}$, it follows that $\bigcup_{h\in \Xi}A_h^c\cap S_{n}$ forms a partition of $S_{n}$. Therefore, we can begin by writing
\[
 |xH\cap S_n| = \sum_{h\in \Xi} |xH\cap S_n\cap A_{h}^c|.
 \]

Using equation (\ref{eq:also_inductive_step}), this becomes
 \[
 |xH\cap S_n| = \sum_{h\in \Xi} |xH\cap S_n\cap A_{h}^c|= \sum_{h\in \Xi} |{{h^{-1}x}}H\cap S_{n-1}\cap A_{h^{-1}}|.
\]
By reindexing, we can replace $h^{-1}$ by $h$ so that we have
 \[
 |xH\cap S_n| = \sum_{h\in \Xi} |{{hx}}H\cap S_{n-1}\cap A_{h}|.
\]

\end{proof}

\begin{lemma}
   Let $H\leq F_2$, $x\in F_2$ and $n\geq 2$. For any $\varnothing\neq \Omega\sub \Xi$,
    \begin{align}\label{inductive_step}
         \sum_{h\in \Omega}|xH\cap S_n\cap A_{h}| =& (|\Omega|-1)\cdot |xH\cap S_n|+\sum_{h\in \Omega^c}|xH\cap S_{n}\cap A_{h}^c |.
    \end{align}

\end{lemma}

\begin{proof}

Recall that $\bigcup_{h\in \Xi}A_h^c$ is a partition of $F_2\setminus \{e\}$, so for each $h\in \Xi$, we have that 
\begin{equation}\label{eq:partition}
A_h = (A_{h}^c)^c = \bigcup_{k\in \Xi\setminus \{h\}}A_k^c.
\end{equation}
When $|\Omega|=1$, the statement follows from equation (\ref{eq:partition}).

Suppose that $|\Omega|>1$. Using equation (\ref{eq:partition}), we have
\begin{align*}
\sum_{h\in \Omega}|xH\cap S_n\cap A_{h}| =& \sum_{h\in \Omega}\left|xH\cap S_n\cap \left(\bigcup_{k\in \Xi\setminus \{h\}}A_{k}^c\right)\right|\\
=&\sum_{h\in \Omega}\left|\bigcup_{k\in \Xi\setminus \{h\}}xH\cap S_n\cap A_{k}^c\right|\\
=&\sum_{h\in \Omega}\sum_{k\in \Xi\setminus \{h\}}\left|xH\cap S_n\cap A_{k}^c\right|\\
=&\sum_{k\in \Xi}\sum_{h\in \Omega\setminus \{k\} }\left|xH\cap S_n\cap A_{k}^c\right|
\end{align*}

where the last equality comes from the fact that
\begin{align*}
&\{(h,k):h\in \Omega, k\in \Xi\setminus \{h\}\} \\=& \{(h,k):h\in \Omega, k\in \Xi, k\neq h\}\\ =&\{ (h,k): k\in \Xi,h\in \Omega\setminus \{k\}\}.
\end{align*}

Now, we will split up the sum over $k\in \Xi$ into a sum over $k\in \Omega$ and a sum over $k\in \Omega^c$.

\begin{align*}
   \sum_{k\in \Xi}\sum_{h\in \Omega\setminus \{k\} }\left|xH\cap S_n\cap A_{k}^c\right|=& \sum_{k\in \Omega}\sum_{h\in \Omega\setminus \{k\} }\left|xH\cap S_n\cap A_{k}^c\right|+\sum_{k\in \Omega^c}\sum_{h\in \Omega\setminus \{k\} }\left|xH\cap S_n\cap A_{k}^c\right|\\
   =& \sum_{k\in \Omega}\sum_{h\in \Omega\setminus \{k\} }\left|xH\cap S_n\cap A_{k}^c\right|+\sum_{k\in \Omega^c}\sum_{h\in \Omega }\left|xH\cap S_n\cap A_{k}^c\right|\\
    =& \sum_{k\in \Omega}(|\Omega|-1)\cdot\left|xH\cap S_n\cap A_{k}^c\right|+\sum_{k\in \Omega^c}|\Omega|\cdot\left|xH\cap S_n\cap A_{k}^c\right|.
    \end{align*}
    We can recombine the above two sums into a sum over $k\in\Xi$ and a sum over $k\in \Omega^c$.
    \begin{align*}
    & \sum_{k\in \Omega}(|\Omega|-1)\cdot\left|xH\cap S_n\cap A_{k}^c\right|+\sum_{k\in \Omega^c}|\Omega|\cdot\left|xH\cap S_n\cap A_{k}^c\right| \\
     =& (|\Omega|-1)\cdot\sum_{k\in \Xi}\left|xH\cap S_n\cap A_{k}^c\right|+\sum_{k\in \Omega^c}\left|xH\cap S_n\cap A_{k}^c\right|\\
    =&(|\Omega|-1)\cdot \left|xH\cap S_n\right|+\sum_{k\in \Omega^c} \left|xH\cap S_n\cap A_{k}^c\right|.
\end{align*}
Notice that the last equality comes from applying equation (\ref{eq:also_inductive_step}) and (\ref{rec_relation}) respectively. Re-indexing the second sum to be in terms of $h$ instead of $k$ gives, as desired,
\begin{equation*}
(|\Omega|-1)\cdot \left|xH\cap S_n\right|+\sum_{h\in \Omega^c} \left|xH\cap S_n\cap A_{h}^c\right|.    
\end{equation*}
 
\end{proof}

Now for the proof of Theorem \ref{thm:main}.

\begin{proof}

Fix $H\leq F_2$, $y\in F_2$ and $n\geq 2$. We will show by induction that there exist constants $a_{i,xH}$ and sets $\Omega_{i,xH}\sub \Xi$, defined for each $i\in \N$ and $xH\in F_2/H$, which depend on $yH$ but not on $n$, such that for each $1\leq k \leq n-1$ we have
\begin{equation}\label{eq:induction_hypothesis}
    |yH\cap S_n|=  \sum_{i=1}^{k}\sum_{xH\in F_2/H} a_{i,xH}\cdot |xH\cap S_{n-i}|+\sum_{\substack{zH\in F_2/H }}\sum_{\substack{h\in \Omega_{k+1,zH}}}|zH\cap S_{n-(k+1)}\cap A_{h} |.
\end{equation}

Setting $k=n-1$ in equation (\ref{eq:induction_hypothesis}) yields the statement of the theorem since $S_0\cap A_h=\varnothing$ for all $h\in\Xi$ as $S_0$ contains only the identity.

We define $\Omega_{-1,xH} = \varnothing$ for all $xH \in F_2/H$ while putting $\Omega_{0,yH} = \Xi$ and $\Omega_{0,xH} = \varnothing$ when $xH\neq yH$. Now for any $xH\in F_2/H$ and $k> 0$, after having defined $\Omega_{k,xH}$ for each $xH\in F_2/H$, we define 
\begin{equation}\label{8}
    \Omega_{k+1,zH} = \{h\in \Xi:\Omega_{k-1,zH}=\varnothing \text{ and }\Omega_{k,h^{-1}zH}\neq \varnothing\}.
\end{equation}

Then for $i\geq 1$ and $xH\in F_2/H$, we define $a_{i,xH} = 0$ when $\Omega_{i,xH}$ is empty and $a_{i,xH} = (|\Omega_{i,xH}|-1)$ when $\Omega_{i,xH}$ is nonempty. 

Continuing, we will show equation (\ref{eq:induction_hypothesis}) for the case $k=0$. From equation (\ref{rec_relation}), we have 
\[
  |yH\cap S_n| = \sum_{h\in \Xi}|hyH\cap S_{n-1}\cap A_{h}|.
\]

Recall that for any $xH\in F_2/H$, we have \begin{align*} 
&\Omega_{1,xH} = \{h\in \Xi:  \Omega_{-1,xH}=\varnothing \text{ and } \Omega_{0,h^{-1}xH}\neq \varnothing\} \\=& 
\{h\in \Xi:\Omega_{0,h^{-1}xH}\neq \varnothing\}= \{h\in \Xi: h^{-1}xH = yH\},
\end{align*}
which in particular means that $h\in \Omega_{1,xH}$ iff $xH = hyH$. Hence, $\{\Omega_{1,xH}\}_{xH\in F_2/H}$ forms a partition of $\Xi$. This gives us that
\begin{align*}
    |yH\cap S_n| = \sum_{h\in \Xi}|hyH\cap S_{n-1}\cap A_{h}| =& \sum_{xH\in F_2/H}\sum_{h\in \Omega_{1,xH}}|hyH\cap S_{n-1}\cap A_{h}|\\
    =& \sum_{xH\in F_2/H}\sum_{h\in \Omega_{1,xH}}|xH\cap S_{n-1}\cap A_{h}|,
\end{align*}

and that is exactly equation (\ref{eq:induction_hypothesis}) for $k=0$.

Now, suppose that equation (\ref{eq:induction_hypothesis}) holds for $k$, that is,
\begin{equation*}
    |yH\cap S_n|=  \sum_{i=1}^{k}\sum_{xH\in F_2/H} a_{i,xH}\cdot |xH\cap S_{n-i}|+\sum_{\substack{zH\in F_2/H}}\sum_{\substack{h\in \Omega_{k+1,zH}}}|zH\cap S_{n-(k+1)}\cap A_{h} |.
\end{equation*}

We will show that equation (\ref{eq:induction_hypothesis}) holds for $k+1$.

When $\Omega_{k+1,zH}$ is nonempty, we will use equation (\ref{inductive_step}) to say that
    \begin{align*}
        &\sum_{i=1}^{k}\sum_{xH\in F_2/H} a_{i,xH}\cdot |xH\cap S_{n-i}|+\sum_{\substack{zH\in F_2/H}}\sum_{\substack{h\in \Omega_{k+1,zH}}}|zH\cap S_{n-(k+1)}\cap A_{h} |\\
        =& 
       \sum_{i=1}^{k}\sum_{xH\in F_2/H} a_{i,xH}\cdot |xH\cap S_{n-i}|\\
       +&\sum_{\substack{zH\in F_2/H\\ \Omega_{k+1,zH}\neq \varnothing}}\left((|\Omega_{k+1,zH}|-1)\cdot |zH\cap S_{n-(k+1)}|+\sum_{h\in \Omega_{k+1,zH}^c}|zH\cap S_{n-(k+1)}\cap A_{h}^c |\right).
    \end{align*}

    Recall that $a_{k+1,zH} = 0$ when $\Omega_{k+1,zH}$ is empty and $a_{k+1,zH} = (|\Omega_{k+1,zH}|-1)$ when $\Omega_{k+1,zH}$ is nonempty, so this becomes
    \begin{align*}
   & 
       \sum_{i=1}^{k}\sum_{xH\in F_2/H} a_{i,xH}\cdot |xH\cap S_{n-i}|\\
       +&\sum_{\substack{zH\in F_2/H \\ \Omega_{k+1,zH}\neq \varnothing}}\left((|\Omega_{k+1,zH}|-1)\cdot |zH\cap S_{n-(k+1)}|+\sum_{h\in \Omega_{k+1,zH}^c}|zH\cap S_{n-(k+1)}\cap A_{h}^c |\right)\\
       =&
    \sum_{i=1}^{k+1}\sum_{xH\in F_2/H} a_{i,xH}\cdot |xH\cap S_{n-i}|
       +\sum_{\substack{zH\in F_2/H \\ \Omega_{k+1,zH}\neq \varnothing}}\sum_{h\in \Omega_{k+1,zH}^c}|zH\cap S_{n-(k+1)}\cap A_{h}^c |.
    \end{align*}

    Now applying equation (\ref{eq:also_inductive_step}), this becomes
    \begin{align*}
    \sum_{i=1}^{k+1}\sum_{xH\in F_2/H} a_{i,xH}\cdot |xH\cap S_{n-i}|
       +\sum_{\substack{zH\in F_2/H \\\Omega_{k+1,zH}\neq \varnothing}}\sum_{h\in \Omega_{k+1,zH}^c}|zH\cap S_{n-(k+1)}\cap A_{h}^c |\\=    \sum_{i=1}^{k+1}\sum_{xH\in F_2/H} a_{i,xH}\cdot |xH\cap S_{n-i}|
       +\sum_{\substack{zH\in F_2/H \\\Omega_{k+1,zH}\neq \varnothing}}\sum_{h\in \Omega_{k+1,zH}^c}|h^{-1}zH\cap S_{n-(k+2)}\cap A_{h^{-1}} |.
    \end{align*}

While putting $wH = h^{-1}zH$ where $wH\in F_2/H$, consider the following reindexing. Starting with
\begin{align*}
&
\{(h,zH):zH\in F_2/H, \Omega_{k+1,zH}\neq \varnothing, h \in \Omega_{k+1,zH}^c\}, \\
\intertext{we rewrite $h \in \Omega_{k+1,zH}^c$ as $h \in \Xi\setminus \Omega_{k+1,zH}$,}
&
\{(h,zH):zH\in F_2/H, \Omega_{k+1,zH}\neq \varnothing, h \notin \Omega_{k+1,zH}\}. \\
\intertext{Next, using equation (\ref{8}), we know if $h\notin\Omega_{k+1,zH}$, then $\Omega_{k-1,zH}\neq\varnothing$ or $\Omega_{k,h^{-1}zH}=\varnothing$. However, we cannot have $\Omega_{k-1,zH}\neq\varnothing$, since this condition, irrespective of $h$ in equation (\ref{8}), forces $\Omega_{k+1,zH}=\varnothing$. Thus, $h\notin\Omega_{k+1,zH}$ tells us we must have $\Omega_{k,h^{-1}zH}=\varnothing$,}
&
\{(h,zH):zH\in F_2/H, \Omega_{k+1,zH}\neq \varnothing, \Omega_{k,h^{-1}zH}=\varnothing\}. \\
\intertext{Recalling we set $wH=h^{-1}zH$, we know $zH=hwH$. This allows us to say}
&
\{(h,hwH):wH\in F_2/H, \Omega_{k+1,hwH}\neq \varnothing, \Omega_{k,wH}=\varnothing\}. \\
\intertext{Using equation (\ref{8}), we find the conditions $\Omega_{k,wH}=\varnothing \text{ and } \Omega_{k+1,hwH}\neq \varnothing$ provide for an $h^{-1}\in \Omega_{k+2,wH}$, so it follows that}
&
\{(h,hwH):wH\in F_2/H, h^{-1}\in \Omega_{k+2,wH}\}.
\end{align*}

Applying this re-indexing, we obtain
    \begin{align*}
    &\sum_{i=1}^{k+1}\sum_{xH\in F_2/H} a_{i,xH}\cdot |xH\cap S_{n-i}|
       +\sum_{\substack{zH\in F_2/H \\\Omega_{k+1,zH}\neq \varnothing}}\sum_{h\in \Omega_{k+1,zH}^c}|h^{-1}zH\cap S_{n-(k+2)}\cap A_{h^{-1}} | \\
       =& \sum_{i=1}^{k+1}\sum_{xH\in F_2/H} a_{i,xH}\cdot |xH\cap S_{n-i}|
       +\sum_{\substack{wH\in F_2/H }}\sum_{h^{-1}\in \Omega_{k+2,wH}}|wH\cap S_{n-(k+2)}\cap A_{h^{-1}} | \\
       =& \sum_{i=1}^{k+1}\sum_{xH\in F_2/H} a_{i,xH}\cdot |xH\cap S_{n-i}|
       +\sum_{\substack{wH\in F_2/H }}\sum_{h\in \Omega_{k+2,wH}}|wH\cap S_{n-(k+2)}\cap A_{h} |,
    \end{align*}

which is equation (\ref{eq:induction_hypothesis}) for $k+1$.

Now setting $k=n-1$ in equation (\ref{eq:induction_hypothesis}),
\begin{align*}
    |yH\cap S_n|=&  \sum_{i=1}^{n-1}\sum_{xH\in F_2/H} a_{i,xH}\cdot |xH\cap S_{n-i}|+\sum_{zH\in F_2/H }\sum_{\substack{h\in \Omega_{n,zH}}}|zH\cap S_{n-(n-1+1)}\cap A_{h} |\\
    =&
     \sum_{i=1}^{n-1}\sum_{xH\in F_2/H} a_{i,xH}\cdot |xH\cap S_{n-i}|+\sum_{zH\in F_2/H }\sum_{\substack{h\in \Omega_{n,zH}}}|zH\cap S_{0}\cap A_{h} |\\
     =& \sum_{i=1}^{n-1}\sum_{xH\in F_2/H} a_{i,xH}\cdot |xH\cap S_{n-i}|
\end{align*}
since $S_0\cap A_{h} = \varnothing$ for all $h\in \Xi$.

\end{proof}

\section{Proof that the Algorithm is Correct}

The proof of Theorem \ref{thm:main} shows not only that a recurrence relation exists, but also gives a method for calculating the coefficients. Namely, the coefficients $(a_{i,xH})$ in Theorem \ref{thm:main} are equal to $|\Omega_{i,xH}|-1$ when $|\Omega_{i,xH}|>0$ and are equal to $0$ otherwise. The sets $(\Omega_{i,xH})$ are defined recursively such that for any $xH$, we have $\Omega_{-1,xH}=\varnothing$ while $\Omega_{0,yH} = \Xi$ and $\Omega_{0,xH} = \varnothing$ for $xH\neq yH$, and 
\begin{equation}\label{eq:rec_relation_for_omega}
    \Omega_{k+1,zH} = \{h\in \Xi:\Omega_{k-1,zH}=\varnothing \text{ and }\Omega_{k,h^{-1}zH}\neq \varnothing\}.
\end{equation}
However, calculating recurrence relations this way is cumbersome and impractical, so we prove that the coefficients can be calculated using the method described in Section \ref{section:introduction} and shown in Section \ref{section:examples}.

With $yH\in F_2/H$ fixed, we define the family of sets $E_i(xH)$ indexed over $i\in\N$ and $xH\in F_2/H$. These sets will contain directed edges in the coset graph of $F_2/H$. Let $E_1$ be the set of all edges of the form $yH\to zH$ where $zH$ is a vertex adjacent to $yH$. Having defined $E_i$, define $E_{i+1}$ by
\begin{align*}
    E_{i+1}= \{tH\to sH: tH &\text{ and } sH \text{ are adjacent vertices}, sH\to tH\not\in E_i, \\
    &\text{ and there exists }rH \text{ adjacent to } tH\text{ such that } rH\to tH\in E_i\}.
\end{align*}

For $xH\in F_2/H$, define $E_i(xH)$ as the set of directed edges that have the form $zH\to xH$ while belonging to $E_i$. These sets $E_i$ describe exactly the directed edges that are highlighted green in Step $i$ as in the examples given in Section \ref{section:examples}.

Now all that is left is to show the following lemma.

\begin{lemma}\label{lem:E_i=omega_i}
    Given any $xH\in F_2/H$, we have $|E_i(xH)| = |\Omega_{i,xH}|$ for all $i\in\N$.
\end{lemma}

\begin{proof}

We will show that the following two statements are true for any $i\in\N$.

\begin{equation}\label{eq:E_in_terms_of_Omega}
E_i(xH) = \{h^{-1}xH\to xH: h\in \Omega_{i,xH}\} \text{ for all } xH\in F_2/H
\end{equation}
which is equivalent to
\begin{equation}\label{equation:Omega_in_terms_of_E}
\Omega_{i,xH} = \{h\in \Xi:  h^{-1}xH\to xH \in E_i(xH)\}\text{ for all } xH\in F_2/H.
\end{equation}

First, let us see why these two equations are equivalent. The coset graph of $F_2/H$ could contain parallel edges, that is, distinct edges which both begin and end at same pair of vertices. This occurs when there is $zH\in F_2/H$ and distinct $h_1,h_2\in \Xi$ such that $h_1^{-1}zH= h_2^{-1}zH$. From the recursive definition of $(E_i)$, it is clear that if two edges are parallel, then one is contained in $E_i(xH)$ iff the other is too. Similarly, if $h_1^{-1}zH= h_2^{-1}zH$ for $h_1,h_2\in \Xi$, then $h_1\in \Omega_{i,zH}$ iff $h_2\in \Omega_{i,zH}$. This shows that equation (\ref{eq:E_in_terms_of_Omega}) is true iff equation (\ref{equation:Omega_in_terms_of_E}) is true, for if either equation is true, we can create a bijection which sends each edge of the form $zH\to xH$ in $E_{i}(xH)$ to an $h\in \Omega_{i,xH}$ with $h^{-1}xH=zH$ or vice versa. Further, this implies that $|\Omega_{i,xH}| = |E_i(xH)|$. Thus, it suffices to prove that equation (\ref{eq:E_in_terms_of_Omega}) holds for all $i\in\N$.

We proceed by induction on $i$. Recall that $\Omega_{-1,zH} = \varnothing$ for all $zH$ while $\Omega_{0,yH}=\Xi$ and $\Omega_{0,zH}=\varnothing $ for $zH\neq yH$. It then follows that, by equation (\ref{eq:rec_relation_for_omega}), we have that $\Omega_{1,xH} = \{h\in \Xi: h^{-1}xH = yH\}$. Given $E_1(xH)$ is the set of all edges of the form $yH\to hyH$ where $hyH=xH$ for some $h\in\Xi$, it directly follows that 
\begin{align*}
E_1(xH) =& \{yH\to xH: xH = hyH\text{ for some } h\in \Xi\} \\
=& \{h^{-1}xH\to xH: h^{-1}xH = yH\text{ for some } h\in \Xi\} \\
=&  \{h^{-1}xH\to xH: h\in \Omega_{1,xH}\}.
\end{align*}
This shows the base case of the induction.

Having fixed $i\in\N$, suppose that 

\begin{equation}
E_i(xH) = \{h^{-1}xH\to xH: h\in \Omega_{i,xH}\} \text{ for all } xH\in F_2/H
\end{equation}
holds, and we will show that 

\begin{equation*}
E_{i+1}(xH) = \{h^{-1}xH\to xH: h\in \Omega_{i+1,xH}\} \text{ for all } xH\in F_2/H.
\end{equation*}

Pick any $xH\in F_2/H$ and rewrite the recursive definition of $E_{i+1}$, so that we have
\begin{align*}
    E_{i+1}(xH)= \{wH\to xH:&~hwH = xH \text{ for some } h\in\Xi \text{ while } xH\to wH\not\in E_i(wH)\\
    &\text{ and } h'wH\to wH\in E_{i}(wH) \text{ for some } h'\in\Xi\}.
\end{align*}

By the induction hypothesis, we know the condition $xH\to wH\not\in E_i(wH)$ is equivalent to $h^{-1}\not\in \Omega_{i,wH}$ whenever $xH=hwH$. Similarly, we have the condition that $h'wH\to wH\in E_{i}(wH)$ for some $h'\in\Xi$ is equivalent to there existing an $h'\in \Xi$ such that $h'^{-1}\in \Omega_{i,wH}$. Thus, it follows that the condition $h'wH\to wH\in E_{i}(wH)$ for some $h'\in\Xi$ is sufficient for $\Omega_{i,wH}\neq\varnothing$, but using the induction hypothesis it becomes clear it is also necessary. Then $hwH=xH$ means that $wH=h^{-1}xH$, and so we can now rewrite the equation for $E_{i+1}(xH)$ as 
\begin{align*}
    E_{i+1}(xH)=\{h^{-1}xH\to xH: h\in\Xi, h^{-1}\not\in\Omega_{i,h^{-1}xH} \text{, and } \Omega_{i, h^{-1}xH}\neq \varnothing\}.
\end{align*}

The condition $ h^{-1}\not\in \Omega_{i,h^{-1}xH}$ means that $\Omega_{i-2,h^{-1}xH}\neq \varnothing$ or $\Omega_{i-1,xH}=\varnothing$. However, the condition $\Omega_{i,h^{-1}xH}\neq\varnothing$ implies that $\Omega_{i-2,h^{-1}xH}=\varnothing$. Together, this means we can rewrite the equation for $E_{i+1}(xH)$ as

\begin{align*}
        E_{i+1}(xH)=&\{h^{-1}xH\to xH: 
(\Omega_{i-2,h^{-1}xH}\neq \varnothing \text{ or } \Omega_{i-1,xH}=\varnothing)\text{ and } 
\Omega_{i,h^{-1}xH}\neq \varnothing\}
    \\=&\{h^{-1}xH\to xH: 
\Omega_{i-1,xH}=\varnothing\text{ and } 
\Omega_{i,h^{-1}xH}\neq \varnothing\}\\
       =&\{h^{-1}xH\to xH: h\in \Omega_{i+1,xH}\}.
\end{align*}
This completes the proof.

\end{proof}

\section{Proof of Theorem \ref{thm:main_2}}\label{section:rec_relation_is_finite}

The goal of this section is to show Theorem \ref{thm:main_2} by showing that the algorithm described in the introduction terminates after finitely many steps. More specifically, we will show the following theorem.

\begin{theorem}\label{thm:rec_relation_is_finite}
    Let $H$ be a finite index subgroup of $F_2$, which contains an element of odd length, and let $y\in F_2$. There is an $N\in \N$ large enough such that the constants $(a_{i,xH})_{i\geq 1, xH\in F_2/H}$ from Theorem \ref{thm:main} satisfy $a_{i,xH}  = 0$ for any $xH$ and any $i\geq N$.
\end{theorem}

Fix a finite index subgroup $H\leq F_2$ and fix $y\in F_2$. Recall that the constants $(a_{i,xH})$ can be derived from the sets $(\Omega_{i,xH})$, which were defined in the proof of Theorem \ref{thm:main} given at the end of Section \ref{section:proof_of_main_thm}. As a consequence of Lemma \ref{lem:E_i=omega_i}, Theorem \ref{thm:rec_relation_is_finite} can be established by showing $E_{i} = \varnothing$ for all large enough $i$, so we will prove this theorem by exploring some properties of the sets $E_i$. To this end, consider the following definition.

\begin{definition}
    Let $(v_0,\dots, v_n)$ be a walk in the coset graph of $F_2/H$. We call this walk \textit{vital} if there does not exist a walk $(u_0,\dots, u_{m})$ in the coset graph of $F_2/H$ satisfying $u_0 = v_0$, $u_{m-1} = v_{n+1}$ and $u_m = v_m$ with $m<n$.
\end{definition}

Below, we will give proofs of the following two facts.

\begin{theorem}\label{thm:vital_means_green}
Let $\ell \in\N$ and let $xH, zH\in F_2/H$. Then the following two statements are true:
\begin{enumerate}[label = (\roman*)]
    \item If $(v_0,\dots, v_{\ell})$ is a vital walk with $v_0=yH$, then $v_{j-1}\to v_{j}\in E_j$ for each $1\leq j\leq \ell$.
    \item If $zH\to xH\in E_{\ell}$, then there is a vital walk $(v_0,\dots, v_{\ell})$ with $v_0 = yH$, $v_{\ell-1}=zH$, and $v_{\ell} = xH$.
\end{enumerate}
\end{theorem}

\begin{theorem}\label{thm:finitely_many_vital}
Suppose that $H$ contains an element of odd length. Then there exists an $N$ such that there are no vital walks of length larger than $N$. 
\end{theorem}

Theorem \ref{thm:rec_relation_is_finite} follows immediately from Theorems \ref{thm:vital_means_green} and \ref{thm:finitely_many_vital} since if $E_i\neq \varnothing$ for infinitely many values of $i$, then Theorem \ref{thm:vital_means_green} gives us arbitrarily long vital walks, which Theorem \ref{thm:finitely_many_vital} says is impossible.

In order to prove Theorem \ref{thm:vital_means_green}, we first need the following lemma.

\begin{lemma}\label{lem:minimial_obstructing_path}
    Suppose that $(v_0,\dots, v_n)$ is a walk that is not vital. Let $m$ be the minimal natural number such that there exists a walk $(u_0,\dots, u_m)$ such that $u_0 = v_0$, $u_{m-1}= v_{m+1}$ and $u_m = v_m$ with $m<n$. Then the walk $(u_0,\dots, u_m)$ is vital.
\end{lemma}

\begin{proof}

Let $(v_0, \ldots, v_n)$ be a walk. Say that a walk $(s_0,\dots, s_j)$ satisfies property $(*)$ if there is a $j<n$ with $s_0 = v_0$, $s_{j-1}= v_{j+1}$ and $s_j = v_j$.

Suppose for the sake of contradiction that $m$ is the minimal natural number such that there exists a walk $(u_0,\dots, u_m)$ that satisfies property $(*)$, but that the walk $(u_0,\dots, u_m)$ is not vital.

Then there is a walk $(w_0,\dots, w_k)$ with $w_0 = u_0$, $w_{k-1} = u_{k+1}$ and $w_{k}=u_k$ with $k<m$. Now, we will create a new walk out of $(w_0,\dots, w_k)$ and $(u_0,\dots, u_m)$, which satisfies property $(*)$ and contradicts the minimality of $m$.

Consider the walk given by $(w_0,\dots, w_{k-2}, u_{k+1},\dots, u_{m})$, which is a walk since $w_{k-2}$ is adjacent to $w_{k-1} = u_{k+1}$. This walk has length $m-1$, starts with $w_0=u_0=v_0$ and has a second to last element $u_{m-1}=v_{m+1}$ as well as a last element $u_m=v_m$.
It follows that this walk satisfies property $(*)$ and so contradicts the minimality of $m$.
Hence it must be that $(u_0,\dots, u_m)$ is vital.
\end{proof}
Now for the proof of Theorem \ref{thm:vital_means_green}.

\begin{proof}

We will proceed by induction on $\ell$. We know that statements (i) and (ii) are both true for $\ell=1$ since $E_1$ consists of all edges pointing out of $yH$ and any walk of length $1$ is vital. Hence it follows that whenever $xH$ is adjacent to $yH$, we have that $yH\to xH$ is contained in $E_1$ and that the walk $(yH,xH)$ is vital.

Now suppose that (i) and (ii) are both true for $\ell$. That is, suppose that if a walk $(v_0,\dots, v_\ell)$ with $v_0 = yH$ is vital, then $v_{j-1}\to v_{j}\in E_j$ for each $1\leq j\leq \ell$ and suppose that if $zH\to xH\in E_{\ell}$, then there is a vital walk $(v_0,\dots, v_{\ell})$ with $v_0 = yH$, $v_{\ell-1} = zH$, and $v_{\ell} = xH$.

We will show that (i) and (ii) are both true for $\ell+1$.

First we will start with (i). Suppose that $(v_0,\dots, v_{\ell+1})$ is a vital walk with $v_0=yH$. From the definition of vital walk, it follows that $(v_0,\dots, v_{\ell})$ is also a vital walk and so the induction hypothesis for (i) gives us that $v_{j-1}\to v_j\in E_j$ for all $1\leq j\leq \ell$. We will show that $v_{\ell}\to v_{\ell+1}\in E_{\ell+1}$. From the recursive definition of $E_{\ell+1}$, we have that $v_{\ell}\to v_{\ell+1}\in E_{\ell+1}$ iff $v_{\ell+1}\to v_{\ell}\not\in E_{\ell}$ and there is also a $rH$ with $rH\to v_{\ell}\in E_{\ell}$. We can take $rH = v_{\ell-1}$ allowing us to say $v_{\ell-1}\to v_{\ell}\in E_{\ell}$. Thus, all that is left to show is that $v_{\ell+1}\to v_{\ell}\not\in E_{\ell}$.

Suppose for the sake of contradiction that $v_{\ell+1}\to v_{\ell}\in E_{\ell}$. By the inductive hypothesis for (ii), this means that there is a vital walk $(u_0,\dots, u_{\ell-1},u_{\ell})$ with $u_0 = yH$, $u_{\ell-1} = v_{\ell+1}$ and $u_{\ell} = v_{\ell}$ that contradicts the fact that $(v_0,\dots, v_{\ell+1})$ is vital. Thus, it must be that $v_{\ell+1}\to v_{\ell}\not\in E_{\ell}$, so we get $v_{\ell}\to v_{\ell+1}\in E_{\ell+1}$, which proves statement (i) for $\ell+1$.

Now for statement (ii), we begin by letting $xH, zH\in F_2/H$. Then we suppose that $zH\to xH\in E_{\ell+1}$, meaning we have that $xH\to zH\not \in E_{\ell}$ and there is an $rH\in F_2/H$ with $rH\to zH\in E_{\ell}$ as well. From the induction hypothesis for (ii), this means that there is a vital walk $(v_0,\dots, v_{\ell})$ with $v_0 = yH$, $v_{\ell-1} = rH$ and $v_{\ell} = zH$. Let $v_{\ell+1} = xH$. We will show that $(v_0,\dots, v_{\ell},v_{\ell+1})$ is a vital walk, which will prove statement (ii).

Suppose for the sake of contradiction that $(v_0,\dots, v_{\ell},v_{\ell+1})$ is not a vital walk. Then there is a walk $(u_0,\dots, u_{m})$ such that $u_0 = yH$, $u_{m-1} = v_{m+1}$ and $u_{m}=v_m$ where $m<\ell+1$. From the fact that $(v_0,\dots, v_{\ell})$ is a vital walk, we know that $m$ cannot be less than or equal to $\ell-1$, which implies $m>\ell-1$, so it must be that $m=\ell$ since $m<\ell+1$. Thus, $m=\ell$ is the only and minimal choice that can still satisfy property (*), and therefore from Lemma \ref{lem:minimial_obstructing_path}, it follows that the walk $(u_0,\dots, u_{\ell})$ is vital. By the inductive hypothesis for (i), this means that $u_{\ell-1}\to u_{\ell}\in E_{\ell}$. Notice that since $u_{\ell-1} = v_{\ell+1} = xH$ and $u_{\ell} = v_{\ell} = zH$, we have $xH\to zH\in E_{\ell}$, but this is a contradiction, for we earlier arrived at $xH\to zH\not\in E_{\ell}$. This means that $(v_0,\dots, v_{\ell+1})$ is a vital walk, and hence statement (ii) is true for $\ell+1$.

By induction, statements (i) and (ii) hold for all $\ell$.

\end{proof}

Now we will set our sights on proving Theorem \ref{thm:finitely_many_vital}, which will follow from this next lemma.

\begin{lemma}\label{lem:erg_thm}
    Let $K$ be any finite index subgroup of $F_2$, and suppose that $K$ contains an element of odd length. Then there is an $N$ large enough such that $K$ contains an element of length $n$ for any $n\geq N$.
\end{lemma}

For this lemma, we require an Ergodic Theorem found in \cite[Theorem 4.10.3, pg. 152]{tempelman}.
\begin{theorem}\label{thm:erg_thm}
    Let $g\mapsto U_g$ be a unitary action of $F_2$ on a Hilbert Space. Let $\Xi$ be the set of generators of $F_2$ along with their inverses. For $n\geq 1$, let $S_n$ be the set of elements of $F_2$ which have length $n$. Pick any $h\in \Xi$. Then for any vector $v$ belonging to the Hilbert Space, we have that 
    \begin{equation*}
        \lim_{n\to\oo}\frac{1}{|S_n\cup hS_n|}\sum_{g\in S_n\cup hS_n}U_gv 
    \end{equation*}
    is the projection of $v$ onto the subspace of vectors which are invariant under each $(U_g)_{g\in F_2}$.
\end{theorem}

\begin{corollary}\label{cor:erg_thn}
    For any subset $A\sub F_2$, any $h\in \Xi$, and for any $n\geq 1$, define 
    \begin{equation*}
        d_n^h(A) = \frac{|A\cap (S_n\cup hS_n)|}{|S_n\cup hS_n|} =  \frac{|A\cap S_n|+|A\cap hS_n|}{2|S_n|}=\frac{1}{2}\left(\frac{|A\cap S_n|}{|S_n|}+\frac{|A\cap hS_n|}{|S_n|}\right).
    \end{equation*}

    Let $K$ be a finite index subgroup of $F_2$. Then for any $h\in \Xi$ and any $z\in F_2$, $\lim_{n\to\oo}d_n^h(zK) = \frac{1}{[F_2:K]}$.
\end{corollary}

\begin{proof}

Consider the vector space of functions from $F_2/K\rightarrow \C$. This is finite dimensional and hence a Hilbert Space. Consider the unitary $F_2$-action given by $g\mapsto U_g$ where $(U_gf)(xK) = f(gxK)$. The subspace of $F_2$-invariant vectors is the space of constant functions and the projection of $f:F_2/K\rightarrow  \C$ onto this subspace is given by $\frac{1}{|F_2/K|}\sum_{xH\in F_2/K}f(xK) = \frac{1}{[F_2:K]}\sum_{xH\in F_2/K}f(xK)$. Then a special case of Theorem \ref{thm:erg_thm} says that for any function $f:F_2/K\rightarrow \C$, we have that
\begin{equation*}
    \lim_{n\to\oo}\frac{1}{|S_n\cup hS_n|}\sum_{g\in S_n\cup hS_n}U_gf(K) = \frac{1}{[F_2:K]}\sum_{xH\in F_2/K}f(xK).
\end{equation*}

Let $z\in F_2$ and let $f = 1_{zK}$. Observe that, for each $n\geq 1$, 
    \begin{equation*}
        \sum_{g\in S_n\cup hS_n}U_gf(K) = \sum_{g\in S_n\cup hS_n}1_{zK}(gK) = |zK\cap (S_n\cup hS_n)| = |zK\cap S_n|+|zK\cap hS_n|.
    \end{equation*}

    Additionally, we have that $|S_n\cup hS_n| = |S_n|+|hS_n| = |S_n|+|S_n| = 2|S_n|$. 
     Using Theorem \ref{thm:erg_thm} we have
\begin{align}\label{eq:erg_thm_eq}
\lim_{n\to\oo}\frac{1}{2}\left(\frac{|zK\cap S_n|}{|S_n|}+\frac{|zK\cap hS_n|}{|S_n|}\right)=&\lim_{n\to\oo}\frac{1}{|S_n\cup hS_n|}\sum_{g\in S_n\cup hS_n}U_gf(K)\\ =& \frac{1}{[F_2:K]}\sum_{xH\in F_2/K}f(xK) = \frac{1}{[F_2:K]}
\end{align}
as desired.

\end{proof}

Now for the proof of Lemma \ref{lem:erg_thm}.

\begin{proof}
We are given at $K$ contains an element of odd length, call it $g$. Write $g = g_{\ell}\cdots g_{1}$, where $g_i\in \Xi$ for all $i$ and $\ell$ is odd. Set $p_0 = e$ and $p_j = g_j\cdots g_1$ for $1\leq j\leq \ell$ so that $g_{j+1}p_j = p_{j+1}$ and $p_{\ell}=g\in K$.

\textbf{Claim:}
$\lim_{n\to\oo}\frac{|p_{j}K\cap S_n|}{|S_n|}$ exists for $0\leq j\leq \ell$.

To prove this claim, first note that the following holds for each $1\leq j\leq \ell-1$, 
\[
\lim_{n\to\oo}d^{g_{j+1}}_n(p_jK) = \lim_{n\to\oo}\frac{\frac{|p_jK\cap S_n|}{|S_n|}+\frac{|g_{j+1}p_jK\cap S_n|}{|S_n|}}{2} = \lim_{n\to\oo}\frac{\frac{|p_jK\cap S_n|}{|S_n|}+\frac{|p_{j+1}K\cap S_n|}{|S_n|}}{2}
\]
Additionally note that $p_{\ell}K = gK = K$.
Then
\begin{align*}
&\sum_{j=0}^{\ell-1}\lim_{n\to\oo}d^{g_{j+1}}_n(p_{j}K) =
\sum_{j=0}^{\ell-1}\lim_{n\to\oo}\frac{\frac{|p_jK\cap S_n|}{|S_n|}+\frac{|p_{j+1}K\cap S_n|}{|S_n|}}{2} 
= \lim_{n\to\oo}\sum_{j=0}^{\ell-1}\frac{\frac{|p_jK\cap S_n|}{|S_n|}+\frac{|p_{j+1}K\cap S_n|}{|S_n|}}{2}. 
\end{align*}
Separating out the first and last terms of this sum gives
\begin{align*}
&\lim_{n\to\oo}\frac{\frac{|K\cap S_n|}{|S_n|}+\frac{|p_1K\cap S_n|}{|S_n|}}{2} + \sum_{j=1}^{\ell-2}\frac{\frac{|p_jK\cap S_n|}{|S_n|}+\frac{|p_{j+1}K\cap S_n|}{|S_n|}}{2}+\frac{\frac{|p_{\ell-1}K\cap S_n|}{|S_n|}+\frac{|K\cap S_n|}{|S_n|}}{2}
\\
=&
\lim_{n\to\oo}\frac{|K\cap S_n|}{|S_n|}+\sum_{j=1}^{\ell-1}\frac{|p_{j}K\cap S_n|}{|S_n|}\\ =& \lim_{n\to\oo}\frac{|K\cap S_n|}{|S_n|}+2\sum_{j=1}^{\frac{\ell-1}{2}}\frac{\frac{|p_{2j-1}K\cap S_n|}{|S_n|}+\frac{|p_{2j}H\cap S_n|}{|S_n|}}{2}.
\end{align*}
However, $\lim_{n\to\oo} 2\sum_{j=1}^{\frac{\ell-1}{2}}\frac{\frac{|p_{2j-1}K\cap S_n|}{|S_n|}+\frac{|p_{2j}K\cap S_n|}{|S_n|}}{2}= 2\sum_{j=1}^{\frac{\ell-1}{2}}\lim_{n\to\oo}d^{g_{2j}}(p_{2j-1}K)$, which implies that $\lim_{n\to\oo}\frac{|K\cap S_n|}{|S_n|}$ exists. The existence of the limit $\lim_{n\to\oo}\frac{|p_{j}K\cap S_n|}{|S_n|}$ implies the existence of the limit $\lim_{n\to\oo}\frac{|p_{j+1}K\cap S_n|}{|S_n|}$ based on the fact that $\lim_{n\to\oo}d^{g_{j+1}}_n(p_jK)$ exists. Hence, by induction, each of the limits  $\lim_{n\to\oo}\frac{|p_{j}K\cap S_n|}{|S_n|}$ exist.

Now one more claim.

\textbf{Claim:}
  $\lim_{n\to\oo}\frac{|K\cap S_n|}{|S_n|} = \frac{1}{[F_2:H]}$

This can be shown by considering 
\[
0 = \lim_{n\to\oo}d_n^{g_1}(K) - d_n^{g_2}(p_1K) = \lim_{n\to\oo}\frac{|K\cap S_n|}{|S_n|} - \lim_{n\to\oo}\frac{|p_2K\cap S_n|}{|S_n|}
\]
so $\lim_{n\to\oo}\frac{|K\cap S_n|}{|S_n|} = \lim_{n\to\oo}\frac{|p_2K\cap S_n|}{|S_n|}$. Likewise from $\lim_{n\to\oo}d^{g_j}(p_{j-1}K) - d_n^{g_{j+1}}(p_jK) = 0$, we eventually find that $\lim_{n\to\oo}\frac{|K\cap S_n|}{|S_n|} = \lim_{n\to\oo}\frac{|p_{\ell-1}K\cap S_n|}{|S_n|}$ since $\ell-1$ is even. Then 
\[
\frac{1}{[F_2:H]}=\lim_{n\to\oo}d^{g_{\ell}}(p_{\ell-1}K) = \lim_{n\to\oo}\frac{\frac{|p_{\ell-1}K\cap S_n|}{|S_n|}+\frac{|K\cap S_n|}{|S_n|}}{2} = \lim_{n\to\oo}\frac{|K\cap S_n|}{|S_n|}.
\]

From this claim it follows that $|K\cap S_n|>0$ for all large enough $n$, which proves the lemma.

\end{proof}

 We can now finally prove Theorem \ref{thm:finitely_many_vital}.

\begin{proof}

From Lemma \ref{lem:erg_thm}, we can pick $N_0$ such that for each $n\geq N_0$, we have $H$ contains an element of length $n$. Let $N = N_0+2\cdot [F_2:H]$. We will show that there are no vital walks of length larger than $N-1$.

Suppose that $(v_0,\dots, v_i)$ is a walk with $v_0 = yH$ and $i\geq N$. We will show that this is not a vital walk by constructing $(u_0,\dots, u_{i-2},u_{i-1})$ with $u_0=v_0$, $u_{i-2} = v_i$ and $u_{i-1}=v_{i-1}$. To this end, we will take a walk from $yH$ to $H$, a walk from $H$ to $H$, a walk from $H$ to $v_{i}$ and the length 1 walk from $v_{i}$ to $v_{i-1}$ appending them all together to create the walk $(u_0,\dots, u_{i-2},u_{i-1})$.

Let $(s_1,\dots, s_{\ell_1})$ be a walk from $yH$ to $H$ and let $(t_1,\dots, t_{\ell_2})$ be a walk from $H$ to $v_{i}$, with $\ell_1,\ell_2\leq [F_2:H]$. As guaranteed by Lemma \ref{lem:erg_thm}, we can find an element of $H$ with length $\ell_3 = i-\ell_1-\ell_2 \geq N-2\cdot [F_2:H]\geq  N_0$. Let this element be written in reduced form as $g_{i-k}g_{i-k-1}\cdots g_{1}$ where $g_{j}\in \Xi$ for all $j$. Then we have a walk from $H$ to $H$ given by $(H, g_1H, g_2g_1H,\dots, g_{\ell_3}g_{\ell_3-1}\cdots g_{1}H)$.

Appending these together, let 
\begin{equation}
    (u_0,\dots, u_{i-2},u_{i-1}) = (s_1,\dots, s_{\ell_1}, g_1H,\dots, g_{\ell_3}g_{\ell_3-1}\cdots g_{1}H, t_2,\dots, t_{\ell_2}, v_{i-1})
\end{equation}
which is a walk since $s_{\ell_1} = H$, $g_{\ell_3}g_{\ell_3-1}\cdots g_{1}H = H = t_1$, and $t_{\ell_2} = v_{i-2}$. The number of vertices in this walk is 
\[
\ell_1+\ell_3+(\ell_2-1)+1 = \ell_1+\ell_2+\ell_3=i.
\]
Additionally, $u_0 = s_1=yH$, $u_{i-2}= t_{\ell_2} = v_i$ and $u_{i-1} = v_{i-1}$. However, this means that the walk $(v_0,\dots, v_i)$ is not vital, completing the proof.

\end{proof}

\section{Further Questions}

Theorem \ref{thm:rec_relation_is_finite} says that when $H\sub F_2$ contains an element of odd length, then there is an $N$ such that the algorithm for computing the recurrence relation for $|yH\cap S_n|$ will end after $N$ steps. However, the proof of Theorem \ref{thm:rec_relation_is_finite} does not give us the ability to find $N$ or easily give bounds on how large $N$ can be. From computing examples, we conjecture that we can always take $N = [F_2:H]$.

\textbf{Question:} Is it true that the algorithm will always end in at most $[F_2:H]$ steps?

It can be shown that $|S_n|$ has size $4\cdot 3^{n-1}$ for $n\geq 1$ \cite[Pg. 150]{tempelman}. Therefore $\frac{|S_{n}|}{|S_{n-i}|} = 3^i$ for all $i\geq 1$. Let $H$ be a finite index subgroup of $F_2$, which contains an element of odd length, and let $y\in F_2$.  From Theorem \ref{thm:main}
\[
  |yH\cap S_n| = \sum_{i=1}^{N}\sum_{xH\in F_2/H}a_{i,xH}\cdot |xH\cap S_{n-i}|
\]
and so
\[
  \frac{|yH\cap S_n|}{|S_n|} = \sum_{i=1}^{N}\sum_{xH\in F_2/H}a_{i,xH}\cdot \frac{|xH\cap S_{n-i}|}{|S_n|} = \sum_{i=1}^{N}\sum_{xH\in F_2/H}\frac{a_{i,xH}}{3^i}\cdot \frac{|xH\cap S_{n-i}|}{|S_{n-i}|}.
  \]
  In light of the second claim in the proof of Lemma \ref{lem:erg_thm}, we have that $\lim_{n\to\oo}\frac{|yH\cap S_n|}{|S_n|} = \frac{1}{[F_2:H]}$ and $\lim_{n\to\oo}\frac{|xH\cap S_{n-i}|}{|S_{n-i}|} = \frac{1}{[F_2:H]} $ for $i\geq 1$, $xH\in F_2/H$. Taking limits, it follows that
   \begin{equation*}
\frac{1}{[F_2:H]} =\sum_{i=1}^{N}\sum_{xH\in F_2/H}\frac{a_{i,xH}}{3^i}\cdot \frac{1}{[F_2:H]}
  \end{equation*}
  and so 
  \begin{equation*}
\sum_{i=1}^{N}\sum_{xH\in F_2/H}\frac{a_{i,xH}}{3^i}=1.
  \end{equation*}
Looking at the examples given in Section 2, we may notice another pattern. In Example 1, the index of the subgroup is $3$ and sum of the coefficients is $1+1+1+1+3 = 7$. In Example 2 the index of the subgroup is 5 and the sum of the coefficients is $1+1+2+2+1+1+3 = 11$. Lastly in Example $3$, the index of the subgroup is $7$ and the sum of the coefficients is $2+2+1+1+3+3+3 = 15$. In each case the coefficients sum to $2\cdot [F_2:H]+1$.

\textbf{Question:} Is it always true that $\sum_{i=1}^{N}\sum_{xH\in F_2/H}a_{i,xH}=2\cdot [F_2:H]+1$?

The second part of Theorem \ref{thm:main} says that if $H$ contains an element of odd length then the algorithm will end after finitely many steps. The condition that $H$ contains an element of odd length is, in general, not necessary. Let $H$ be the subgroup consisting of all elements with even length, $H = \langle a^2, b^2, ab\rangle$. $H$ has index $2$ and it can be verified that the algorithm ends after just one step, $|H\cap S_n| = 3|aH\cap S_{n-1}|$. 

This makes it natural to ask if the odd length condition can be removed entirely.

\textbf{Question}: Is there a subgroup $H\sub F_2$ such that $[F_2:H]<\oo$, but the algorithm does not end after finitely many steps?

\end{document}